\documentclass[11pt,a4paper]{article}
\usepackage[T1]{fontenc}
\usepackage[utf8]{inputenc}
\usepackage{lmodern,microtype,amsmath,amssymb,graphicx,float,caption}
\usepackage[margin=28mm]{geometry}
\usepackage[hidelinks]{hyperref}
\title{A quantitative two-interval criterion\\for the Riemann hypothesis}
\author{Payam Danesh}
\date{}
\begin{document}
\maketitle
\begin{abstract}
We aim to study whether the Weil form can detect a zero of the Riemann zeta function away from the critical line using two functions with increasingly narrow, widely separated supports. We prove that the Riemann hypothesis is equivalent to the Cauchy-Schwarz inequality on every such pair at all sufficiently large separations, provided the component widths decrease subexponentially. A zero to the right of the critical line forces the normalized cross term to exceed an exponential at arbitrarily large separations, with any exponent smaller than the horizontal displacement of that zero. In particular, reciprocal-width intervals suffice. Our proof combines a short-support positivity estimate, a smooth partition of fixed test functions and a pole obstruction for a one-sided Laplace transform. We also express the cross term through prime powers in a moving multiplicative interval. The resulting criterion specifies a uniform arithmetic inequality whose proof would establish the Riemann hypothesis.
\end{abstract}

Keywords: Riemann hypothesis, Weil quadratic form, explicit formula, localized test functions, Laplace transform.

\section{Introduction}

In November 1859, Bernhard Riemann in his paper states that every nontrivial zero of the zeta function has real part \(1/2\). Weil's explicit formula connects these zeros with prime powers and leads to a positivity criterion on smooth test functions \cite{ref1}. The question addressed here is how severely the geometry of the test supports can be restricted while retaining a criterion equivalent to the Riemann hypothesis (RH). The analytic ingredients include the classical zero-counting estimate and the functional equation, as developed in Titchmarsh's account \cite{ref2}.

Through the previous works, Yoshida studied Hermitian forms attached to zeta functions and established positivity on sufficiently short intervals \cite{ref3}. Bombieri also developed a variational analysis of the Weil form, including functions supported on finite unions of intervals \cite{ref4}. These results supply an existing setting for disconnected supports and local positivity. Connes and Consani studied archimedean positivity through the semilocal trace formula and showed that finitely many transform-vanishing conditions can be imposed in Weil's criterion \cite{ref5}. Suzuki subsequently related the growth of a function built from zeta zeros to the analytic continuation of its one-sided transform \cite{ref6}. The pole argument below belongs to this analytic tradition and Its role is to recover information about zeros from a bound imposed only on pairs of functions with shrinking supports. In a related spectral formulation, a study associates a quadratic energy to each zeta zero and express its imaginary part through the horizontal displacement from the critical line \cite{ref7}. Here that displacement controls a lower exponential growth rate for a normalized Weil cross term. In another work, Bombieri and Lagarias \cite{ref8} used the Guinand--Weil explicit formula to connect Li's positivity criterion with that of Weil.

We focus on the normalized Weil cross term for two test functions, each with positive Weil energy, supported on intervals of prescribed width that shrink as the intervals move farther apart. The main question is whether functions supported in \(\left\lbrack 0,w(L) \right\rbrack\) and \(\left\lbrack L - w(L),L \right\rbrack\) suffice, where \(L\) is the outer span and \(w(L)\) decreases to zero. Further, we try to prove that they do whenever \(\log\left( 1/w(L) \right) = o(L)\). A zero with real part \(\beta > 1/2\) forces violations larger than \(e^{\eta L}\) at arbitrarily large spans for each \(0 < \eta < \beta - 1/2\). The reciprocal-width case \(w(L) = 1/L\) follows directly.

\section{The Weil form and positivity}

To start, let us \(\zeta\) denote the meromorphic continuation of \(\sum_{n \geq 1}^{}n^{- s}\) from \(\mathfrak{R}s > 1\) and let \(\mathcal{D} =C_{c}^{\infty}\left( \mathbb{R};\mathbb{C} \right)\) denote the complex-valued smooth functions of compact support, with \(L^{2}\) norm \(\Vert f \Vert_{2} = \left( \int_{\mathbb{R}}^{}\left| f(x) \right|^{2}\, dx \right)^{1/2}\). The support is the closure of the set on which a function is nonzero. We let \(Z\) be the multiset of nontrivial zeros \(\rho = \beta + i\gamma\) of \(\zeta\), with multiplicities, so \(0 < \beta < 1\). The centered zero and bilateral Laplace transform can be defined by

\begin{equation}
z_{\rho} = \rho - \frac{1}{2},\quad\quad H_{f}(z) = \int_{\mathbb{R}}^{}f(x)e^{zx}\, dx\quad\left( f\in \mathcal{D},\ z\in \mathbb{C} \right).
\label{eq:1}
\end{equation}

All coordinates and parameters are dimensionless here. The centered Weil form, linear in its first argument, is

\begin{equation}
B(f,g) = \sum_{\rho \in Z}^{}H_{f}\left( z_{\rho} \right)\overline{H_{g}\left( - \overline{z_{\rho}} \right)},\quad\quad Q(f) = B(f,f).
\label{eq:2}
\end{equation}

For each nonnegative integer \(r\), integration by parts gives

\begin{equation}
\left| H_{f}(\sigma + iy) \right| \leq C_{f,r}\left( 1 + |y| \right)^{- r},\quad\quad|\sigma| \leq \frac{1}{2},
\label{eq:3}
\end{equation}

uniformly in \(\sigma,y\). The zero-counting estimate

\begin{equation}
N(T): = \#\{\rho \in Z:\left| \mathfrak{I}\rho \right| \leq T\} \ll (1 + T)\log(2 + T)
\label{eq:4}
\end{equation}

holds for \(T \geq 0\) and counts both signs of the ordinates \cite{ref2}. Together these bounds prove absolute convergence in (2). Here \(U \ll V\), or \(U = O(V)\), means \(|U| \leq CV\) in the stated range and subscripts record allowed dependence of \(C\). For positive \(U,V\), the notation \(U \asymp V\) means that both \(U \ll V\) and \(V \ll U\) hold. Reflection under \(z \mapsto - \bar{z}\) and reindexing show that \(B(f,g) = \overline{B(g,f)}\), so \(Q(f)\) is real. For the translation \(\tau_{a}f(x) = f(x - a)\), \(a\in \mathbb{R}\), we have

\begin{equation}
H_{\tau_{a}f}(z) = e^{az}H_{f}(z),\quad\quad B\left( \tau_{a}f,\tau_{a}g \right) = B(f,g).
\label{eq:5}
\end{equation}

The exponential factors cancel exactly in the second identity.

We prove a short-support lower bound to justify the denominator used below. The unitary Fourier transform and the archimedean multiplier are

\begin{equation}
\widehat{f}(y) = \frac{1}{\sqrt{2\pi}}\int_{\mathbb{R}}^{}f(x)e^{- iyx}\, dx,\quad\quad\Omega(y)= \mathfrak{R}\frac{\Gamma'}{\Gamma}\left( \frac{1}{4} + \frac{iy}{2} \right) - log\pi,
\label{eq:6}
\end{equation}

which \(\Gamma\) is the gamma function. We let \(k_{f}(u) = \int_{\mathbb{R}}^{}f(x + u)\overline{f(x)}\, dx\). The von Mangoldt function is \(\Lambda\left( p^{j} \right) = logp\) for primes \(p\) and positive integers \(j\) and zero on all other positive integers. The complete explicit formula is

\begin{equation}
\begin{aligned}
Q(f) &= 2\mathfrak{R}\left( H_{f}\left( \frac{1}{2} \right)\overline{H_{f}\left( - \frac{1}{2} \right)} \right) + \int_{\mathbb{R}}^{}\Omega(y)\left| \widehat{f}(y) \right|^{2}\, dy \\
 & - 2\mathfrak{R}\sum_{n \geq 2}^{}\frac{\Lambda(n)}{\sqrt{n}}\, k_{f}\left( \log n \right).
\end{aligned}
\label{eq:7}
\end{equation}

This is the additive Guinand--Weil formula \cite{ref4} which the first term comes from the poles and the prime sum is finite by compact support. Stirling's estimate and continuity of \(\Omega\) give \(C_{0} > 0\) such that

\begin{equation}
\Omega(y) \geq log\left( e + |y| \right) - C_{0}\quad\quad\left( y\in \mathbb{R} \right).
\label{eq:8}
\end{equation}

\textbf{Lemma 1.} If \(f\in \mathcal{D}\) has support in an interval of length \(0 < d < log2\), then

\begin{equation}
Q(f) \geq \left\lbrack \frac{1}{2}\log\left( e + \frac{\pi}{2d} \right) - C_{0} - 2de^{d/2} \right\rbrack \Vert f \Vert_{2}^{2}.
\label{eq:9}
\end{equation}

There is consequently \(d_{0} > 0\) such that \(Q(f) > 0\) for every nonzero \(f\in \mathcal{D}\) supported in an interval of length at most \(d_{0}\).

\emph{Proof.} we can translate the support into \(\lbrack - d/2,d/2\rbrack\) using Eq (5). Cauchy--Schwarz gives \(\left| H_{f}( \pm 1/2) \right| \leq d^{1/2}e^{d/4} \Vert f \Vert_{2}\), so the pole term in Eq. (7) is at least \(- 2de^{d/2} \Vert f \Vert_{2}^{2}\). Since \(k_{f}(u) = 0\) for \(|u| > d\), all prime terms vanish. The pointwise estimate \(\left| \widehat{f}(y) \right| \leq \left( d/(2\pi) \right)^{1/2} \Vert f \Vert_{2}\) implies, for \(M = \pi/(2d)\),

\begin{equation}
\int_{- M}^{M}\left| \widehat{f}(y) \right|^{2}\, dy \leq \frac{dM}{\pi} \Vert f \Vert_{2}^{2} = \frac{1}{2} \Vert f \Vert_{2}^{2}.
\label{eq:10}
\end{equation}

By Plancherel, at least half the squared norm lies outside this interval. Equations~(8) and~(10) therefore give us

\begin{equation}
\int_{\mathbb{R}}^{}\Omega(y)\left| \widehat{f}(y) \right|^{2}\, dy \geq \left\lbrack \frac{1}{2}\log(e + M) - C_{0} \right\rbrack \Vert f \Vert_{2}^{2}.
\label{eq:11}
\end{equation}

Adding the pole estimate proves Eq (9). Its bracket tends to positive infinity as \(d \downarrow 0\), proving the final assertion. We fix \(L_{*} > 0\) and a nonincreasing function \(w:\left\lbrack L_{*},\infty) \rightarrow (0,1 \right\rbrack\) such that

\begin{equation}
w(L) \rightarrow 0,\quad\quad\log\left( 1/w(L) \right) = o(L)\quad(L \rightarrow \infty).
\label{eq:12}
\end{equation}

The latter condition means that \(\log\left( 1/w(L) \right)/L \rightarrow 0\) and no continuity of \(w\) is assumed. We define

\begin{equation}
I_{L,w} = \left\lbrack 0,w(L) \right\rbrack,\quad\quad J_{L,w} = \left\lbrack L - w(L),L \right\rbrack.
\label{eq:13}
\end{equation}

\textbf{Figure~1} shows the reciprocal-width geometry. The outer span determines the separation, whereas each diagonal energy is evaluated on one short interval.

\begin{figure}[H]
\centering
\includegraphics[width=.92\linewidth]{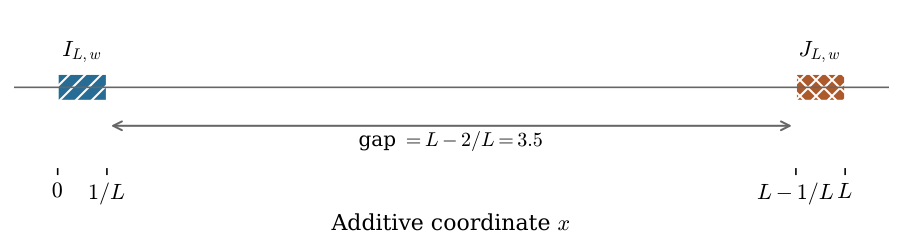}
\caption{Prescribed support intervals at \(L = 4\) for \(w(L) = 1/L\), with hatching distinguishing the intervals of length \(\tfrac{1}{4}\).}
\end{figure}

For sufficiently large \(L\), these intervals are disjoint and \(w(L) \leq d_{0}\). A pair \((f,g)\) is admissible at \(L\) if \(f,g\in \mathcal{D}\) are nonzero and supported in \(I_{L,w},J_{L,w}\), correspondingly. Lemma~1 makes

\begin{equation}
R(f,g) = \frac{\left| B(f,g) \right|}{\sqrt{Q(f)Q(g)}}
\label{eq:14}
\end{equation}

well defined on such pairs.

\section{Localization of a translated test}

In this section, the next estimate transfers a bound on admissible pairs to two translates of a fixed function. The number of pieces depends on the width and keeping this dependence explicit allows the sub exponential extension.

\textbf{Lemma 2.} We assume Eq (12) and let \(\eta \geq 0\) and suppose \(R(f,g) \leq e^{\eta L}\) for every real admissible pair at all sufficiently large \(L\). If \(h\in \mathcal{D}\) is real and \(supph \subset \lbrack - A,A\rbrack\) with \(A \geq 1\), we place

\begin{equation}
C_{A} = 2A + 2,\quad\quad m(t) = \left\lceil \frac{4}{w\left( t + C_{A} \right)} \right\rceil.
\label{eq:15}
\end{equation}

For all sufficiently large \(t\),

\begin{equation}
\left| B\left( \tau_{t}h,h \right) \right| \leq C_{h,A,\eta}\, m(t)\log\left( 2 + m(t) \right)e^{\eta t},
\label{eq:16}
\end{equation}

with a constant independent of \(t\). In particular, for each \(\varepsilon > 0\),

\begin{equation}
\left| B\left( \tau_{t}h,h \right) \right| = O_{h,A,\eta,\varepsilon,w}\left( e^{(\eta + \varepsilon)t} \right).
\label{eq:17}
\end{equation}

For reciprocal width, Eq (16) gives \(\left| B\left( \tau_{t}h,h \right) \right| \ll_{h,A,\eta}t\log(2 + t)e^{\eta t}\).

\emph{Proof.} we choose a real nonnegative bump \(\theta \in C_{c}^{\infty}\left( ( - 1,1) \right)\) positive on \(\lbrack - 1/2,1/2\rbrack\). Its periodic sum \(S(v) = \sum_{j\in \mathbb{Z}}^{}\theta(v - j)\) is smooth and positive. Then \(\chi = \theta/S\) is supported in \(( - 1,1)\) and \(\sum_{j}^{}\chi(v - j) = 1\). For fixed \(t\) , we write \(m = m(t)\) and decompose \(h\) as

\begin{equation}
u_{j}(x) = h(x)\chi(mx - j),\quad h = \sum_{j\in \mathbb{Z}}^{}u_{j},\quad suppu_{j} \subset \left\lbrack (j - 1)/m,(j + 1)/m \right\rbrack.
\label{eq:18}
\end{equation}

Only \(|j| \leq Am + 1\) can contribute, giving \(O_{A}(m)\) nonzero pieces. We convert each piece by \(- j/m\). By substituting \(v = mx\) in its transform gives, uniformly for \(|\sigma| \leq 1/2\), we have

\begin{equation}
H_{\tau_{- j/m}u_{j}}(\sigma + iy) = \frac{1}{m}\int_{- 1}^{1}h\left( (v + j)/m \right)\chi(v)e^{\sigma v/m}e^{iyv/m}\, dv.
\label{eq:19}
\end{equation}

The amplitude and its first two derivatives in \(v\) are bounded independently of \(m,j,\sigma\). The direct bound for \(|y| \leq m\) and two integrations by parts for \(|y| > m\) give us

\begin{equation}
\left| H_{\tau_{- j/m}u_{j}}(\sigma + iy) \right| \leq \frac{C_{h}}{m}\left( 1 + |y|/m \right)^{- 2}.
\label{eq:20}
\end{equation}

The boundary terms vanish since \(\chi\) is smooth and compactly supported in \(( - 1,1)\).

By applying Eq (20) to both transform factors in Eq (2) and using Eq (5), we obtain

\begin{equation}
\left| Q\left( u_{j} \right) \right| \leq \frac{C_{h}'}{m^{2}}\sum_{\rho \in Z}^{}\left( 1 + \left| \mathfrak{I}\rho \right|/m \right)^{- 4} \ll_{h}\frac{\log(2 + m)}{m}.
\label{eq:21}
\end{equation}

For completeness, (4) bounds the contribution below height \(m\) to the sum by \(O\left( m\log(2 + m) \right)\). The zeros with \(2^{r}m < \left| \mathfrak{I}\rho \right| \leq 2^{r + 1}m\), \(r \geq 0\), contribute at most \(O\left( 2^{- 3r}m\log\left( 2 + 2^{r + 1}m \right) \right)\). Summing over these integer \(r\) proves the last bound. We take \(u_{j}\) in the left copy and \(\tau_{t}u_{k}\) in the right. After translating the left endpoint \((j - 1)/m\) to zero, their outer span is

\begin{equation}
L_{jk} = t + (k - j + 2)/m,\quad\quad t - C_{A} \leq L_{jk} \leq t + C_{A}.
\label{eq:22}
\end{equation}

The bounds follow from \(|j|,|k| \leq Am + 1\) and \(m \geq 4\). For large \(t\), every \(L_{jk}\) exceeds the threshold in the hypothesis. Since \(w\) is nonincreasing, Eq (15) gives

\begin{equation}
2/m \leq \frac{1}{2}w\left( t + C_{A} \right) \leq w\left( L_{jk} \right).
\label{eq:23}
\end{equation}

Thus each translated nonzero pair is admissible at its own span \(L_{jk}\). Lemma~1 ensures positive diagonal energies for large \(t\). Hermitian symmetry and the hypothesis actually imply

\begin{equation}
\left| B\left( \tau_{t}u_{k},u_{j} \right) \right| \leq e^{\eta L_{jk}}\sqrt{Q\left( u_{j} \right)Q\left( u_{k} \right)} \ll_{h,A,\eta}e^{\eta t}\frac{\log(2 + m)}{m}.
\label{eq:24}
\end{equation}

Zero pieces contribute nothing and summing Eq (24) over the \(O_{A}\left( m^{2} \right)\) ordered pairs proves Eq (16). Finally, \(4/w\left( t + C_{A} \right) \leq m(t) \leq 5/w\left( t + C_{A} \right)\) because \(w \leq 1\). Eq~(12) gives \(\log m(t) = o(t)\), and hence

\begin{equation}
\log\left( m(t)\log\left( 2 + m(t) \right) \right) = o(t).
\label{eq:25}
\end{equation}

This proves Eq (17). For \(w(L) = 1/L\), the definition of \(m(t)\) gives \(m(t) \asymp_{A}t\), proving the specialization.

The width condition includes \(L^{- a}\) for every fixed \(a > 0\) and \(\exp\left( - L^{b} \right)\) for every fixed \(0 < b < 1\), on a sufficiently large initial domain. In \textbf{Figure~2}, three examples of these analytic functions is compared.

\begin{figure}[H]
\centering
\includegraphics[width=.90\linewidth]{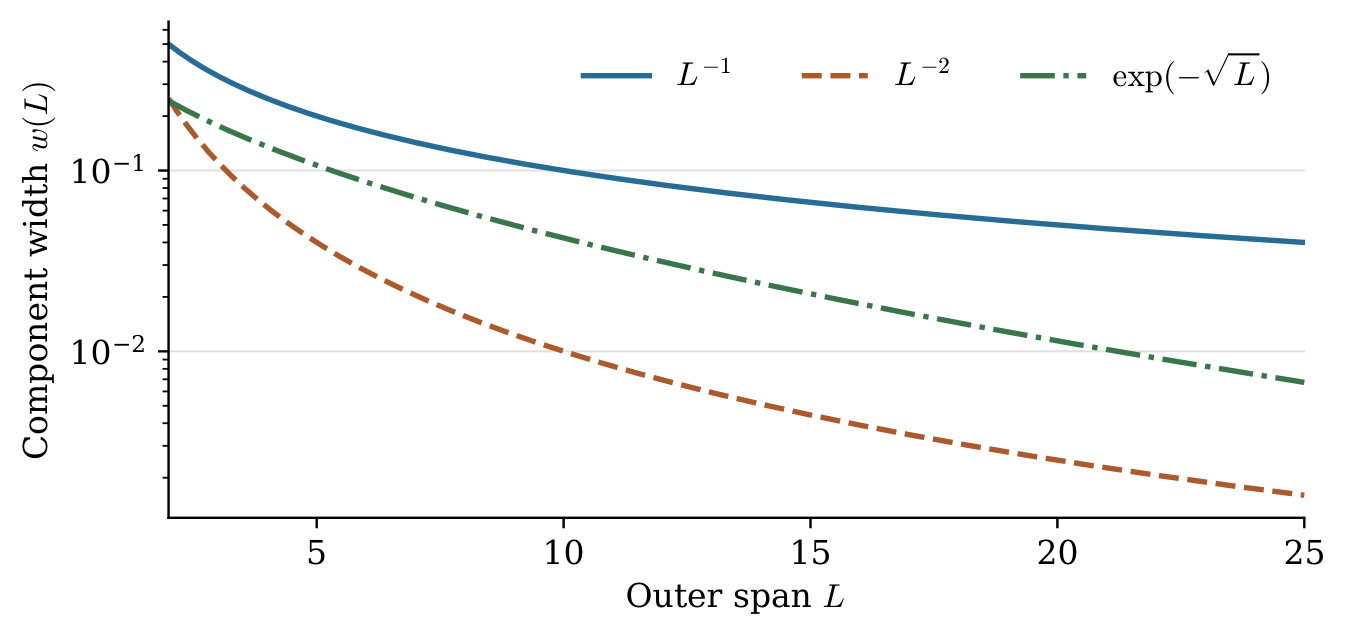}
\caption{Analytic component widths \(L^{- 1}\), \(L^{- 2}\) and \(\exp\left( - \sqrt{L} \right)\) for \(2 \leq L \leq 25\), plotted on a logarithmic vertical axis to illustrate three instances of Eq (12).}
\end{figure}

\section{The criterion and quantitative violations}

A translated even test gives a Laplace transform whose poles retain the centered zeta zeros. This supplies the obstruction to the bound in Lemma~2.

\textbf{Lemma 3.} We suppose \(z_{0} = \rho_{0} - 1/2\) satisfies \(\mathfrak{R}z_{0} > \alpha \geq 0\). There is a real even \(h\in \mathcal{D}\) such that \(F_{h}(t): = B\left( \tau_{t}h,h \right)\) fails to be \(O\left( e^{\alpha t} \right)\) as \(t \rightarrow \infty\).

\emph{Proof.} For real even \(h\), the identities \(H_{h}( - z) = H_{h}(z)\) and \(\overline{H_{h}\left( \bar{z} \right)} = H_{h}(z)\) reduce Eq (2) to

\begin{equation}
F_{h}(t) = \sum_{\rho \in Z}^{}H_{h}\left( z_{\rho} \right)^{2}e^{tz_{\rho}}.
\label{eq:26}
\end{equation}

The coefficient is an ordinary complex square. So we choose a real even bump \(b\in \mathcal{D}\) of integral one and set \(h(x) = r^{- 1}b(x/r)\). Then \(H_{h}\left( z_{0} \right) = H_{b}\left( rz_{0} \right) \rightarrow 1\) as \(r \downarrow 0\). We fix a positive \(r\) with \(H_{h}\left( z_{0} \right) \neq 0\). The series in Eq (26) converges uniformly on bounded \(t\)-intervals. For \(\mathfrak{R}s > 1/2\), equations~(3)-(4) imply

\begin{equation}
\sum_{\rho \in Z}^{}{\int_{0}^{\infty}\left| H_{h}\left( z_{\rho} \right) \right|^{2}}e^{- \left( \mathfrak{R}s- \mathfrak{R}z_{\rho} \right)t}\, dt \leq \frac{\sum_{\rho \in Z}^{}\left| H_{h}\left( z_{\rho} \right) \right|^{2}}{\mathfrak{R}s - 1/2} < \infty.
\label{eq:27}
\end{equation}

Absolute integrability therefore justifies term wise integration and gives us

\begin{equation}
\int_{0}^{\infty}F_{h}(t)e^{- st}\, dt = \sum_{\rho \in Z}^{}\frac{H_{h}\left( z_{\rho} \right)^{2}}{s - z_{\rho}}.
\label{eq:28}
\end{equation}

On each compact, we set disjoint from the centered zeros, the tail with large \(\left| \mathfrak{I}\rho \right|\) converges uniformly by Eq (3). Only finitely many other terms remain. Thus the right side is meromorphic on \(\mathbb{C}\), with residue

\begin{equation}
m_{0}H_{h}\left( z_{0} \right)^{2} \neq 0
\label{eq:29}
\end{equation}

at \(z_{0}\), where \(m_{0}\) is the multiplicity of \(\rho_{0}\). Distinct zeros have distinct pole locations, even when they share an ordinate, therefore they cannot cancel this residue.

If \(F_{h}(t) = O\left( e^{\alpha t} \right)\), the left side of Eq (28) would be holomorphic for \(\mathfrak{R}s > \alpha\). For each compact subset and each derivative order \(n\), its differentiated integrand is bounded for large \(t\) by a constant times \(t^{n}e^{- \epsilon t}\) with some \(\epsilon > 0\). The identity theorem for meromorphic functions extends Eq (28) from \(\mathfrak{R}s > 1/2\) throughout \(\mathfrak{R}s > \alpha\) and this contradicts Eq (29), because \(\mathfrak{R}z_{0} > \alpha\).

\textbf{Theorem 4.} We fix a width function satisfying Eq (12);

\begin{enumerate}
\def\labelenumi{\arabic{enumi}.}
\item
  If the Riemann hypothesis (RH) holds, then \(R(f,g) \leq 1\) for every admissible complex pair at every sufficiently large \(L\).
\item
  If \(\rho_{0} = \beta + i\gamma\) is a nontrivial zero with \(\beta > 1/2\), then for every \(0 < \eta < \beta - 1/2\) and every \(L_{0}\) there exist \(L > L_{0}\) and a real admissible pair at \(L\) for which
\end{enumerate}

\begin{equation}
R(f,g) > e^{\eta L}.
\label{eq:30}
\end{equation}

Consequently the RH is equivalent to \(R(f,g) \leq 1\) for all real admissible pairs at all sufficiently large \(L\).

\emph{Proof.} Under the RH, \(- \overline{z_{\rho}} = z_{\rho}\). Eq~(2) becomes the inner product of the square-summable sequences \(\left( H_{f}\left( z_{\rho} \right) \right)_{\rho \in Z}\) and \(\left( H_{g}\left( z_{\rho} \right) \right)_{\rho \in Z}\). Cauchy-Schwarz proves the first assertion. Here we fix \(\eta\) in the range of the second assertion. If Eq (30) failed for every real pair above some \(L_{0}\), Lemma~2 would apply. We choose \(\varepsilon > 0\) with \(\eta + \varepsilon < \beta - 1/2\). Equation~(17) would then give us \(F_{h}(t) = O\left( e^{(\eta + \varepsilon)t} \right)\) for every real \(h\in \mathcal{D}\), contradicting Lemma~3 with \(\alpha = \eta + \varepsilon\). If the RH fails, reflection across the critical line gives a zero with \(\beta > 1/2\). The second assertion rules out the eventual bound \(R \leq 1\), proving the equivalence.

For the reciprocal-width case, we set \(w(L) = 1/L\). The theorem then applies to \(\lbrack 0,1/L\rbrack\) and \(\lbrack L - 1/L,L\rbrack\) with the same exponent range. The width assumption enters only through the partition cost Eq (25). For \(w(L) = e^{- cL}\) with fixed \(c > 0\), that cost has an exponential rate, so this argument does not give the full exponent range above. A violating pair also produces a negative linear combination of its two components.

\textbf{Corollary 5.} For a real admissible pair satisfying Eq (30), set \(u = f/\sqrt{Q(f)}\) and \(v = g/\sqrt{Q(g)}\). There is \(\epsilon \in \{ - 1,1\}\) such that

\begin{equation}
Q(u - \epsilon v) = 2\left( 1 - R(f,g) \right) < 2\left( 1 - e^{\eta L} \right) < 0.
\label{eq:31}
\end{equation}

\emph{Proof.} Conjugating and reindexing Eq (2) by the conjugation symmetry of \(Z\) shows that \(B(f,g)\) is real for real \(f,g\). We choose \(\epsilon\) to be the sign of \(B(u,v)\), which is nonzero. Expanding the form and using \(Q(u) = Q(v) = 1\) and \(\left| B(u,v) \right| = R(f,g)\) proves Eq (31). This normalization fixes the two individual Weil energies and supplies no bound on \(\Vert u - \epsilon v \Vert_{2}\), so comparison with a negative Rayleigh quotient requires an additional norm estimate.

\hypertarget{arithmetic-formulation}{%
\section{\texorpdfstring{5 Arithmetic formulation }{5 Arithmetic formulation }}\label{arithmetic-formulation}}

For an admissible pair, we define the cross-correlation

\begin{equation}
k(u) = \int_{\mathbb{R}}^{}g(x + u)\overline{f(x)}\, dx,\quad\quad suppk \subset \left\lbrack L - 2w(L),L \right\rbrack.
\label{eq:32}
\end{equation}

It is smooth with compact support and its bilateral Laplace transform is \(H_{g}(z)\overline{H_{f}\left( - \bar{z} \right)}\). For large \(L\) it vanishes on \(( - \infty,0\rbrack\). To derive the arithmetic expression, we put \(\phi(x) = x^{- 1/2}k\left( \log x \right)\) for \(x > 0\) and its multiplicative Mellin transform is

\begin{equation}
\mathcal{M}\phi(s) = \int_{0}^{\infty}\phi(x)x^{s - 1}\, dx = \int_{\mathbb{R}}^{}k(u)e^{(s - 1/2)u}\, du.
\label{eq:33}
\end{equation}

The Guinand--Weil formula, with \(\phi\) zero on \((0,1\rbrack\), reduces to

\begin{equation}
\sum_{\rho \in Z}^{}\mathcal{M}\phi(\rho)= \mathcal{M}\phi(0)+ \mathcal{M}\phi(1) - \sum_{n \geq 2}^{}\Lambda(n)\phi(n) - \int_{1}^{\infty}\phi(x)\frac{x}{x^{2} - 1}\, dx.
\label{eq:34}
\end{equation}

The terms involving \(\phi(1)\) and \(\phi(1/x)\) vanish. In the last integral, \(x = e^{u}\) gives us

\begin{equation}
\int_{1}^{\infty}\phi(x)\frac{x}{x^{2} - 1}\, dx = \int_{0}^{\infty}k(u)\frac{e^{- u/2}}{1 - e^{- 2u}}\, du.
\label{eq:35}
\end{equation}

The support stays away from zero, so no singularity requires regularization. Combining Eqs (32)--(35) yields

\begin{equation}
\begin{aligned}
B(g,f) &= \int_{L - 2w(L)}^{L}k(u)\left( e^{u/2} + e^{- u/2} - \frac{e^{- u/2}}{1 - e^{- 2u}} \right)du \\
 & - \sum_{e^{L - 2w(L)} \leq n \leq e^{L}}^{}\frac{\Lambda(n)}{\sqrt{n}}\, k\left( \log n \right).
\end{aligned}
\label{eq:36}
\end{equation}

Thus the cross term uses a moving interval of prime powers. Each diagonal term has no prime shifts once \(w(L) < log2\), while retaining the pole and archimedean terms in (7).

Theorem~4 identifies the bound needed on the complete right side of Eq (36), uniformly over all admissible pairs. The prime sum and pole integral must be controlled together at the accuracy set by \(\sqrt{Q(f)Q(g)}\). The existence proof also leaves the first violating separation and an algorithm for constructing witnesses from a hypothetical zero undetermined. The result in this work establishes that two subexponentially shrinking supports retain the full logical strength of Weil positivity and quantifies the growth that an off-line zero forces in their normalized interaction. A proof of the Riemann hypothesis through this criterion requires the uniform arithmetic estimate described above.

\section*{Dedication}

I dedicate this work to the memory of Professor Bernhard Riemann. When I visited Biganzolo, I took the photo below. What stayed with me was the last line on his memorial stone: ``\emph{For those who love God, all things must work together for good.} '' Standing before those words I thought of the man whose ideas changed how we understand the prime numbers. I deeply offer this paper in gratitude for the many hours I have spent with his work and for the pleasure I still find in returning to it.

\begin{figure}[H]
\centering
\includegraphics[width=2.817in]{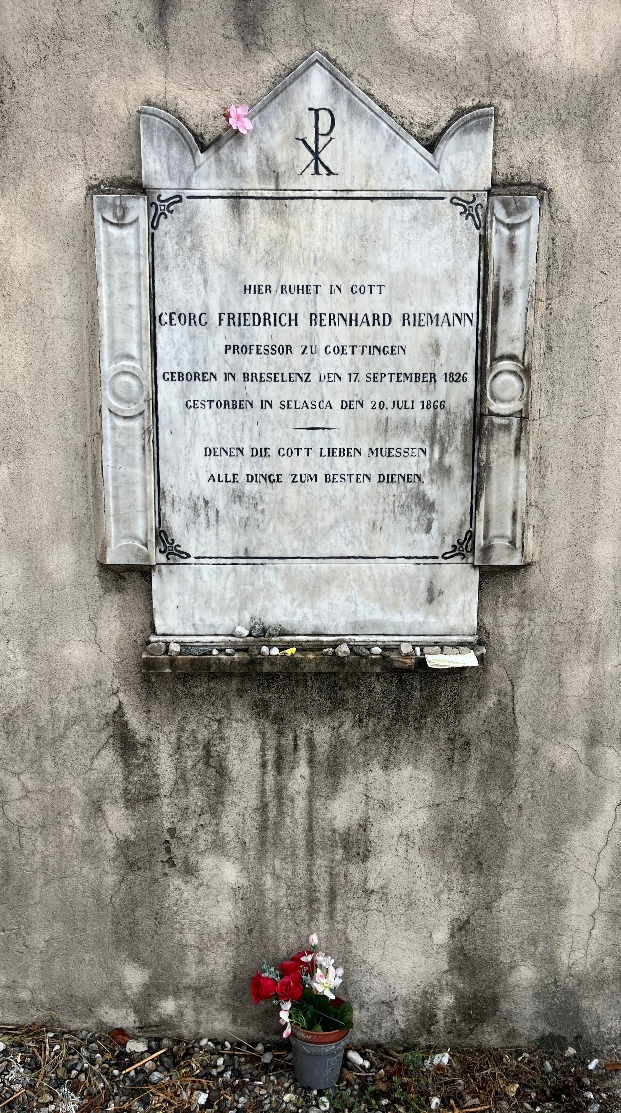}
\caption*{\emph{Bernhard Riemann's memorial stone at Biganzolo Cemetery, Verbania, Italy. 2 September 2026}}
\end{figure}

Polytechnic of Turin

Payam.danesh@polito.it

\begin{thebibliography}{8}

\bibitem{ref1} A. Weil, ``SUR LES FORMULES EXPLICITES DE LA THÉORIE DES NOMBRES,'' \emph{Mathematics of the USSR-Izvestiya}, vol. 6, no. 1, p. 1, Feb. 1972, doi: \nolinkurl{10.1070/IM1972V006N01ABEH001866}.

\bibitem{ref2} E. C. Titchmarsh and D. R. Heath-Brown, ``THE THEORY OF THE RIEMANN ZETA-FUNCTION SECOND EDITION REVISED BY CLARENDON PRESS· OXFORD 1986''.

\bibitem{ref3} H. Yoshida, ``On Hermitian Forms attached to Zeta Functions,'' \emph{https://doi.org/10.2969/aspm/02110281}, vol. 21, pp. 281--326, Jan. 1992, doi: \nolinkurl{10.2969/ASPM/02110281}.

\bibitem{ref4} Bombieri, Enrico. "Remarks on Weil's quadratic functional in the theory of prime numbers, I." \emph{Atti della Accademia Nazionale dei Lincei. Classe di Scienze Fisiche, Matematiche e Naturali. Rendiconti Lincei. Matematica e Applicazioni} 11.3 (2000): 183-233.

\bibitem{ref5} A. Connes and C. Consani, ``Weil positivity and trace formula the archimedean place,'' \emph{Selecta Mathematica 2021 27:4}, vol. 27, no. 4, pp. 77-, Jul. 2021, doi: \nolinkurl{10.1007/S00029-021-00689-4}.

\bibitem{ref6} M. Suzuki, ``Aspects of the screw function corresponding to the Riemann zeta-function,'' \emph{J. Lond. Math. Soc.}, vol. 108, no. 4, pp. 1448--1487, Oct. 2023, doi: \nolinkurl{10.1112/JLMS.12785;JOURNAL:JOURNAL:14697750;WGROUP:STRING:PUBLICATION}.

\bibitem{ref7} P. Danesh, ``Riemann Zeros as quantum energy levels,'' Jul. 2026, doi: \nolinkurl{10.33774/COE-2026-SCGGH}.

\bibitem{ref8} Bombieri, E., \& Lagarias, J. C. (1999). Complements to Li's Criterion for the Riemann Hypothesis. \emph{Journal of Number Theory}, \emph{77}(2), 274--287. \url{https://doi.org/10.1006/JNTH.1999.2392}.
\end{thebibliography}
\end{document}